\documentstyle[amstex,12pt]{amsart}

\renewcommand{\baselinestretch}{1.3}

\makeatletter
\newcommand{\single}{\let\CS=\@currsize\renewcommand{\baselinestretch}{1}\tiny\CS}
\newcommand{\singles}{\let\CS=\@currsize\renewcommand{\baselinestretch}{1.2}\tiny\CS}
\newcommand{\oneanda}{\let\CS=\@currzsize\renewcommand{\baselinestretch}{1.2}\tiny\CS}
\newcommand{\doubles}{\let\CS=\@currsize\renewcommand{\baselinestretch}{1.4}\tiny\CS}

\newcommand{\tree}{\let\CS=\@currsize\renewcommand{\baselinestretch}{1.5}\tiny\CS}
\newcommand{\four}{\let\CS=\@currsize\renewcommand{\baselinestretch}{2}\tiny\CS}

\begin{document} 
\title
[Multiplicity of the multi-graded extended Rees algebra]
{\large {\bf  
A  FORMULA FOR THE MULTIPLICITY OF THE MULTI-GRADED REES ALGEBRA}}

\author{  Clare D'Cruz}
%\\
% Chennai Mathematical Institute,
%92- G. N. Chetty Road,
%T. Nagar, Chennai
%600 017,  INDIA}
\date{}

\newcommand{\se}{\setcounter{equation}{0}}
\newcommand{\ncom}{\newcommand}
\ncom{\bq}{\begin{equation}}
\ncom{\eq}{\end{equation}}
\ncom{\beqn}{\begin{eqnarray*}}
\ncom{\eeqn}{\end{eqnarray*}}
\ncom{\beq}{\begin{eqnarray}}
\ncom{\eeq}{\end{eqnarray}}
\ncom{\been}{\begin{enumerate}}
\ncom{\eeen}{\end{enumerate}}
\ncom{\nno}{\nonumber}
\ncom{\hs}{\mbox{\hspace{.25cm}}}
\ncom{\rar}{\rightarrow}
\ncom{\lrar}{\longrightarrow}
\ncom{\Rar}{\Rightarrow}
\ncom{\noin}{\noindent}
\newtheorem{thm}{Theorem}[section]
\newtheorem{lemma}[thm]{Lemma}
\newtheorem{cor}[thm]{Corollary}
\newtheorem{pro}[thm]{Proposition}
\newtheorem{example}[thm]{Example}
\newtheorem{remark}[thm]{Remark}
\newtheorem{definition}[thm]{Definition}
\ncom{\bt}{\begin{thm}}
\ncom{\et}{\end{thm}}
\ncom{\bl}{\begin{lemma}}
\ncom{\el}{\end{lemma}}
\ncom{\bco}{\begin{cor}}
\ncom{\eco}{\end{cor}}
\ncom{\bp}{\begin{pro}}
\ncom{\ep}{\end{pro}}
\ncom{\bex}{\begin{example}}
\ncom{\eex}{\end{example}}
\ncom{\brm}{\begin{remark}}
\ncom{\erm}{\end{remark}}
\ncom{\bd}{\begin{definition}}
\ncom{\ed}{\end{definition}}
\ncom{\bc}{\begin{center}}
\ncom{\ec}{\end{center}}

\ncom{\comx}{{\mathbb{C}}}
\ncom{\ze}{{\mathbb{Z}}}
\ncom{\re}{{\mathbb{R}}}
\ncom{\Q}{{\mathbb{Q}}}
\ncom{\N}{{\mathbb{N}}}

\ncom{\sz}{\scriptsize}
\ncom{\CM}{Cohen-Macaulay }
\ncom{\sop}{system of parameters}
\ncom{\eop}{\hfill{$\Box$}}
\ncom{\tfae}{the following are equivalent:}
\ncom{\f}{\frac}
\ncom{\la}{\lambda}
\ncom{\si}{\sigma}
\ncom{\ssize}{\scriptsize}
\ncom{\al}{\alpha}
\ncom{\be}{\beta}
\ncom{\Si}{\Sigma}
\ncom{\ga}{\gamma}
\ncom{\kbar}{\overline{\kappa}}
\ncom{\bib}{\bibitem}
\ncom{\sst}{\subset}
\ncom{\sms}{\setminus}
\ncom{\seq}{\subseteq}
\ncom{\est}{\emptyset}
\ncom{\bighs}{\hspace{.5 cm}}
\ncom{\ulin}{\underline}
\ncom{\olin}{\overline}
\ncom{\bip}{\bigoplus}
\ncom{\sta}{\stackrel}

\def\b{{\cal B}}
\def\c{{\cal C}}
\def\r{{\cal R}}
\def\n{{\mathcal N}}
\def\mm{{\cal M}}
\def\ii{{ I_1}}
\def\jj{{ J_1}}
\def\i{{\bf I}}
\def\t{{\bf t}}
\def\m{{\frak m}}
\def\II{{\mathcal I}}

\maketitle
\vspace{-.5in}
\begin{center}
 Chennai Mathematical Institute,\\
92- G. N. Chetty Road,\\
T. Nagar, Chennai 600 017,
India
\end{center}
\begin{center}
{\small email:clare@cmi.ac.in}
\end{center}

\section{\bf INTRODUCTION}

Our aim in this paper is to obtain a formula for the multiplicity of
the maximal homogeneous ideal of the multi-graded extended Rees
algebra. This formula generalizes the 
one obtained in \cite{katz-ver} by Katz and Verma for the ordinary
extended Rees algebra.

The extended Rees algebra was introduced by Rees in \cite{rees-1} and
has been of interest in the past decade.  The Rees algebra which is a
subring of  the extended Rees algebra have been extensively studied. 
It is natural to expect these two algebras to
share some  ring-theoretic properties. For example,
 if $I$ is an ideal of  positive height  in a
Cohen-Macaulay  local ring $R$, then extended Rees algebra of $I$ is
Cohen-Macaulay whenever the Rees algebra is \cite{huneke}. 

Recently, the multi-graded Rees algebra has been  investigated
 (\cite{hyry}, \cite{hhrt}, \cite{ribbe}, \cite{verma},
\cite{clare2}). In particular, the multiplicity of the maximal
homogeneous ideal of the multi-graded Rees algebra was obtained
independently  by Verma (\cite[Theorem~1.4]{verma})
and Herrmann~ et.~al.  (\cite[Corollary~4.7]{hhrt}).
Since the  multi-graded Rees algebra is a subring of the multi-graded
extended Rees algebra, it is natural to study the  multi-graded
extended Rees algebra. In this paper we concentrate on the
multiplicity of the maximal homogeneous ideal of this ring.

 Throughout this paper  $(R, \m)$ will denote a Noetherian local ring of
positive dimension $d$ with infinite residue field.  Let  $I_{1},
\ldots, I_{g}$ be ideals of positive height in $R$ and let $t_{1},
\ldots, t_{g}$ be indeterminates.   The {\em multi-graded extended
Rees algebra} of $R$ with respect to $I_{1}, \ldots,  I_{g}$ is the
graded ring  $\bip_{(r_{1}, \ldots, r_{g}) \in \ze^g}  
   (I_{1} t_{1})^{r_{1}} \cdots (I_{g} t_{g})^{r_{g}}$ 
and will  be denoted by  ${\b}({\bf I})$. 
Here $I_i^{r_i} = R$ if $r_i \leq 0$.  
 Let
${\n}({\i})
= (I_{1} t_{1}, \ldots, I_{g} t_{g}, \m, t_{1}^{-1},
   \ldots, t_{g}^{-1})$.  
The
{\em multi-graded Rees algebra} of $R$ with respect to $I_1, \ldots,
I_g$ 
is the graded ring 
$\bip_{(r_{1}, \ldots, r_{g}) \in  \N^g} 
   (I_{1} t_{1})^{r_{1}} \cdots (I_{g} t_{g})^{r_{g}}$ 
and will  be denoted by  ${\r}({\i})$.
 Let
$ \mm(\i)$ be the maximal homogeneous ideal of $\r({\i})$.
When $g=1$, we say $\b(I)$ is the extended Rees algebra and $\r(I)$ is
the Rees algebra. 

To state our main result we need to define mixed multiplicities. 
Let $I_{1}$ be an $\m$-primary ideal and let $I_{2}, \ldots, I_{g}$ be
ideals of positive height in $(R,\m)$. Then 
for $r_{1}, \ldots, r_{g}$ large, 
$\ell_R (I_{1}^{r_{1}}   I_{2}^{r_{2}} \cdots I_{g}^{r_{g}}/
         I_{1}^{r_{1}+1} I_{2}^{r_{2}} \cdots I_{g}^{r_{g}})$ 
is a polynomial of degree $d-1$ and the terms of degree $d-1$ are
captured in the following sum
$$
  \sum_{ q_{1} + \cdots + q_{g} = d-1}
 e(I_{1}^{[q_{1}+1]} | I_{2}^{[q_{2}]} |\cdots | I_{g}^{[q_{g}]})
  {r_{1} + q_{1} \choose q_{1}} \cdots {r_{g} + q_{g} \choose q_{g}}.
$$
Here 
$e(I_{1}^{[q_{1}+1]} | \cdots | I_{g}^{[q_{g}]})$ 
are positive integers and they are called the {\em mixed
multiplicities} of the set of ideals ($I_{1}, \ldots, I_{g}$)
\cite{teis}.

We state the multiplicity formula for the multi-graded Rees algebra.
In this paper, if $(R,\m)$ is a local ring, then $e(I)$ will denote
the multiplicity of an $\m$-primary ideal in $R$ and $e(R)$ will
denote the multiplicity of the maximal ideal $\m$ of $R$. 
\vspace{.2in}

\bt
{\em (}\cite[Theorem~1.4]{verma}, \cite[Corollary~4.7]{hhrt} {\em )}  
Let $I_{1}, \ldots, I_{g} \seq \m$ be ideals of positive height in
$R$. Then
$$
e(\r({\i})_{\mm(\i)})
= \sum_{q + q_{1} + \cdots +  q_{g} = d-1} 
  e(\m^{[q+1]} | I_{1}^{[q_{1}]} | \cdots | I_{g}^{[q_{g}]}).
$$
\et

The main result of this paper is:

\vspace{.2in}

\bt
\label{mult}
Let $I_{1}, \ldots, I_{g} \seq \m$ be ideals of positive height in $R$. Put 
$L=\m^2 + I_{1} + \cdots + I_{g}$. Then
\beqn
   e( \b( {\bf I})_{{\n}({\i})})
= \f{1}{2^d} \left[ 
  \sum_{t=0}^{g} 
  \sum_{\sta{\scriptstyle q + q_{1}+ \cdots + q_t = d-1}
       {1 \leq i_{1} <  \cdots < i_t \leq g}}
   2^{q_{1} +  \cdots +  q_{t}}
   e({L^{[q+1]} | I_{i_{1}}^{[q_{1}]} |  \cdots | 
                  I_{i_t}^{[q_{t}]}}) \right]  .
\eeqn
\et  

To attain our goal, we need to express mixed multiplicities of certain
homogeneous ideals in the extended Rees ring in terms of mixed
multiplicities of ideals in the ring $R$ (Proposition~\ref{l4}).   We
 recover the multiplicity formula obtained by Katz and Verma  for the
extended Rees algebra (see
Corollary~\ref{extra2}).

We now describe the organization of this paper. In Section three we
prove our main result. Section two is devoted to develop the necessary
preliminary results.   We end this paper by explicitly stating the
multiplicity formula for 
${\b}(\i)_{{\n}({\i})}$ when  $d = 1,2$.

\singles

\vspace{.2in}

\section{  \bf PRELIMINARIES}
In this section we   prove  a few basic facts.

\vspace{.2in}

\bl
\label{bin}
Let $r, s$ and $n$ be non-negative integers. Then 
$$
  \sum_{i=0}^{n} {i+r \choose r} {n-i+s \choose s} 
= {n+r+s+1 \choose r+s+1}.
$$
\el
\pf Apply induction on $n+r$. \qed

\vspace{.2in}

\bl
\label{prod}
Let $s$ be a positive integer and let $r$ and $n$ be  indeterminates. Then 
\beqn
{nr+s \choose s} = n^s {r+s \choose s} 
                  + f_{1}(n) {r+s-1 \choose s-1} 
                  + \cdots + f_{s}(n) 
\eeqn
where $f_{1}(n), \ldots, f_{s}(n) \in \Q(n)$.
\el
\pf Since $\{ {r+i \choose i} \}_{i \in \N}$ form a basis of $\Q(n)[r]$
as a vector space over $\Q(n)$, we can write
\beqn
     {nr+s \choose s} 
&=& \sum_{i=0}^{s} f_{i}(n)  {r+s-i \choose s-i}
\eeqn
where $f_{i}(n) \in \Q(n)$, $0 \leq i \leq s$. Comparing the
coefficient of $r^s$ we get $n^s = f_0(n)$. \qed

\vspace{.2in}
                 
\bl
\label{kmv}
Let $K$ be 
an $\m$-primary ideal and let $J_{1}, \ldots, J_{g} \seq \m$
be ideals of positive height in $R$. Then for  $r_{1},
\ldots, r_{g}$  large,  
$$
  F(r_{1}, \ldots, r_{g})
= \ell \left( \f{J_{1}^{r_{1}} \ldots J_{g}^{r_{g}}}
               {K J_{1}^{r_{1}} \ldots J_{g}^{r_{g}}} \right)
$$
is a polynomial of total degree at most $d-1$ in $r_{1}, \ldots, r_{g}$.
\el
\pf For all large values of $r_{1}, \ldots, r_{g}$, $F(r_{1}, \ldots, r_{g})$
is a polynomial say $P(r_{1}, \ldots, r_{g}) \in \Q[r_{1}, \ldots, r_{g}]$
\cite[Theorem~4.1]{hhrt}.  Since the monomials of highest degree in
$r_{1},  \ldots, r_{g}$  have non-negative coefficients, the total degree
of $P(r_{1}, \ldots, r_{g})$ is equal to the degree of $P(r, \ldots, r)$
and  
\beqn
\begin{array}{llll}
&&       \mbox{deg}~P(r, \ldots, r)  &\\
&=& \dim~ {\displaystyle 
           \left(\f{R[(J_{1} \cdots J_{g}) t]}
                {K R[(J_{1} \cdots J_{g}) t]}
              \right)} - 1
       \hspace{.5in}  & \mbox{[by \cite[Theorem~4.1, Lemma~1.1]{hhrt}]}\\
&\leq& \dim~R[(J_{1} \cdots J_{g}) t] - 2 
                      & \mbox{[by \cite[Theorem~1.5]{rees0}]}\\
&=& d-1 
        &       \mbox{[by \cite[Corollary~1.6]{valla}]}.
\end{array}
\eeqn 
\qed

\vspace{.2in}

\bl
\label{nl2}
Let $J$ be an $\m$-primary ideal and let $I \seq \m$ be an ideal of
positive height in $R$. Put $T = \b(I)$, $M = (t^{-1}, J,
It)$,  $K = J^2 + I$ and $H= J+I$. Then for $0 \leq j < r$,  
\beqn
&&   M^{2(r-j)} \\
&=& \bip_{i=0}^{\infty} R t^{-(2r+i)}
    \bip \bip_{i=1}^{r-1} K^{i-j} t^{-(2r-2i)}
    \bip \bip_{i=0}^{j-1} Rt^{-(2r-2i-1)}
    \bip \bip_{i=j}^{r-1} H K^{i-j} t^{-(2r-2i-1)}  \bip K^{r-j}\\
&&  \bip \bip_{i=0}^{\infty}  (It)^{(2r+i)}
    \bip \bip_{i=1}^{r-1} K^{i-j} (It)^{(2r-2i)}
    \bip \bip_{i=0}^{j-1} (It)^{2r-2i-1}
    \bip \bip_{i=j}^{r-1} H K^{i-j} (It)^{2r-2i-1}.
\eeqn
\el
\pf  
An expression for $M^n$ has been obtained in
\cite[Lemma~3.1]{katz-ver}. By arranging and re-indexing the terms we
get the expression in the above form. 
\qed

\vspace{.2in}

\section{ \bf THE MAIN THEOREM}
In this section we prove our main result.  One of the main ingredients
is Proposition~\ref{l4}.   

\vspace{.2in}

\noin
{\bf Notation:}
Put  ${\b}_0 :=R$, ${\n}_0 := \m$, 
$ L = \m^2 + I_{1} + \cdots + I_{g}$ and  $I_{g+1} = (0)$.  
For $j=1, \ldots, g$ we  inductively define 
${\b}_{j} = {\b}_{j-1}[I_{j} t_j, t_j^{-1}]$, 
${\n}_{j} = (t_{j}^{-1}, {\n}_{j-1}, I_{j} {\b}_{j-1}t_{j})$ and 
$L_{j} = \n_{j}^2 + (I_{j+1} + \cdots + I_{g}) \b_{j}$.

\vspace{.2in}

\bl
\label{dim1}
Let $(R,\m)$ be a local ring of positive dimension $d$. Let $I_1,
\ldots I_g$ be ideals of positive height in $R$. Then $\dim~ \b({\i}) =
\dim~R + g$. 
\el
\pf Since ${\b}_{g} = {\b}_{g-1}[I_{g} t_g, t_g^{-1}]$, it is enough to
show that $\dim~ {\b}_{1} =\dim~R + 1$ and that the ideal
$I_2{\b}_{1}$ has positive height.  
Clearly   ${\b}_{1}/(t_1^{-1})$ is the associated graded ring $G = \oplus_{n \geq
0} I_1^n / I_1 ^{n+1}$ and  $\dim~G = \dim R$. Hence,   $\dim~ \b(I) =  \dim~G + 1 =
\dim~R +1$. 
If the height of $I_2{\b}_{1}$ is zero, then it is contained in some
minimal prime of ${\b}_1$. By a result of Valla,
\cite[cf. Proposition 1.1 (iii)]{valla}, 
$I_2 = I_2{\b}_{1} \cap R$ is contained in some
minimal prime of $R$  which leads to a contradiction. \qed

\vspace{.2in}

\bl
\label{dim2}
Let $I$ be an ideal of positive height in  $(R, \m)$. 
 Let ${\II}_1$  be an $\n(I)$-primary homogeneous ideal and let
 ${\II}_2, \ldots, {\II}_g$  be
homogeneous ideals of positive height in $\b(I)$.
Then  for all non-negative integers $q_1, \ldots, q_g$ with 
$q_1 + \cdots + q_g = d-1$ we have
\beqn
 e(({\II_1}_{\n(I)})^{[q_1 + 1]}| ({\II_2}_{\n(I)})^{[q_2]} 
      | \ldots |  ({\II_g}_{\n(I)})^{[q_g]} )
= e(\II_1^{[q_1 + 1]} | \II_2^{[q_2]} | \ldots |\II_g^{[q_g]}).
\eeqn
\el
\pf Since $\n(I)$ is the maximal homogeneous ideal of $\b(I)$  and
 ${\II}_2, \ldots, {\II}_g$  are homogeneous ideals the following
 isomorphism holds true: 
\beqn
      \f{\II_1^{r_1}     \II_2^{r_2} \ldots \II_g^{r_g}}
        {\II_1^{r_1 + 1} \II_2^{r_2} \ldots \II_g^{r_g}}
\cong \f{(\II_1^{r_1} \II_2^{r_2} \ldots \II_g^{r_g})_{\n(I)}}
        {(\II_1^{r_1+1} \II_2^{r_2} \ldots \II_g^{r_g})_{\n(I)}}.
\eeqn
\qed

\vspace{.2in}

\noin
{\bf Notation:}
Let $I_1$ be an $\m$-primary ideal and let $I_2$ be an ideal of
positive height in $(R, \m)$. 
For  $g=2$, we will use the notation 
\beqn
e_{q}(I_{1} | I_{2}) &:=& e(I_{1}^{[d-q]} | I_{2}^{[q]})
\hspace{2in} q=0, \ldots , d-1,  \\
e_{d}(I_{1} | I_{2}) &:=& 0. 
\eeqn

\vspace{.2in}

\bp
\label{l4}
Let $J$ be an $\m$-primary ideal and  let $I,\ii \seq \m$ be ideals of
positive height in $R$. Let $\jj \seq J+I$ be any ideal of $R$.  Put $T =
\b(I)$, $M = (t^{-1}, J, It)$ and  $K = J^2 + I$. Then for all
$q=0, \ldots, d$ 
\beqn
  e_{q}( M^2 + \jj T|\ii T) 
=2 \left[ e_{q}( K + \jj|\ii)
 + \sum_{q_0+q_{1}+q=d-1} 2^{q_{1}}
        e((K + \jj)^{[q_0+1]}|I^{[q_{1}]}|\ii^{[q]}) \right]. 
\eeqn
\ep

\vspace{.2in}

\vspace{.2in}

The following lemma is well known and easy to see, but nevertheless we
mention  it for the sake of completion.

\vspace{.2in}

\bl
\label{dim3}
Let $I$ be an $\m$-primary ideal and let $J$ be any ideal of positive
height in a local ring of positive dimension $d$. Then for all $q=0,
\ldots d-1$, 
\been 
\item
\label{one}
${\displaystyle
e_{q}(I^{r} | J^{s}) = r^{d -q} s^q e_q(I|J)}.$

\item
$e_0(I|J) = e(I)$. 
\eeen
\el
\pf It is easy to see that  for $r,s \gg 0$
\beqn
\ell \left( 
\f{ (I^n)^r  (J^m)^s}{ (I^n)^{r+1} (J^m)^s} \right)
= \sum_{i=0}^{n-1}
\ell \left( 
\f{I^{nr + i}  J^{ms}}{ (I^{nr+ i+1} J^{ms}} \right).
\eeqn
Comparing the coefficient of $r^{d-1-q} s^{q}$ on both sides we get
(\ref{one}). 
 The second result was proved by D.~Katz and J.~ Verma in
 \cite[Lemma~2.2]{katz-ver}.  
\qed

\vspace{.2in}

\bco
\label{extra1}
Let $J$ be an $\m$-primary ideal and  let $I,\ii \seq \m$ be ideals of
positive height in $R$. Put $T = \b(I)$, 
$M = (t^{-1}, J, It)$. Then for all
$q=0, \ldots, d$ 
\beqn
     e_{q} (M|\ii T)
&=& \f{1}{2^{d-q}}
    \left[ e_{q} (J^2 + I | \ii)
  + \sum_{q_0 + q_{1} + q = d-1} 2^{q_{1}}
     e((J^2 + I)^{[q_0 + 1]} | I^{[q_{1}]} | \ii^{[q]}) \right].
\eeqn
\eco
\pf Since $\dim~\b(I) = d+1$,  
$
e_{q} (M^2| \ii)  = 2^{d+1 - q} e_{q}(M | \ii T)
$ 
by Lemma~\ref{dim3}. 
 Put $\jj= (0)$ in Proposition~\ref{l4}. Then for all
$q=0, \ldots, d$ 
\beqn
e_{q} (M^2| \ii)
= 2  \left[ e_{q} (J^2 + I | \ii)
  + \sum_{q_0 + q_{1} + q = d-1} 2^{q_{1}}
     e((J^2 + I)^{[q_0 + 1]} | I^{[q_{1}]} | \ii^{[q]}) \right].
\eeqn
\qed

\vspace{.2in}

\noin
{\bf Proof of Proposition~\ref{l4}:}  
First note that by Lemma~\ref{dim2}, 
\beqn
\f{((M^2 + \jj T)^r (\ii T)^s)_{\n(I)}}
  {((M^2 + \jj T)^{r+1} (\ii T)^s)_{\n(I)}}
= \f{(M^2 + \jj T)^r (\ii T)^s}
{(M^2 + \jj T)^{r+1} (\ii T)^s}.
\eeqn
Since  $\dim~ \b(I) = d+1$,  for $r,s \gg 0$,  
$\ell(((M^2 + \jj T)^r(\ii T)^s)_{\n(I)}
     /((M^2 + \jj T)^{r+1}(\ii T)^s)_{\n(I)})$ is a 
polynomial of total degree $d$ in $r$ and $s$ \cite{bhat} and can be
written in the form
\beq
\label{ns5}
  \ell \left( 
  \f{(M^2 + \jj T)^r (\ii T)^s_{\n(I)}}
    {(M^2 + \jj T)^{r+1} (\ii T)^s_{\n(I)}} \right) 
= \sum_{q=0}^{d} e_{q}(M^2 + \jj T| \ii T) 
  {r + d-q \choose d-q}{s+q \choose q} + \cdots.
\eeq
But $(M^2 + \jj T)^r \ii T^s$ is a graded ideal. 
and  for all non-negative  integers $r$ and $s$,  the module
$(M^2 + \jj T)^r(\ii T)^s/(M^2 + \jj T)^{r+1}(\ii T)^s$ can
 be expressed as a (finite)  direct sum of $R$- modules which 
have finite length. 
 
Put $H= I+J$. 
 Notice that $(\ii T)^s = \ii^s T$.  It follows from Lemma~\ref{nl2}
 that 
\beqn
&&   (M^2 + \jj T)^r (\ii T)^s\\
&=&  \bip_{i=0}^{\infty} \ii^s t^{-(2r+i)}
     \bip \bip_{i=1}^{r-1} \ii^s (K + \jj)^i t^{-(2r-2i)} 
     \bip \bip_{i=0}^{r-1} \ii^s H (K + \jj)^i  t^{-(2r-2i-1)} \bip \ii^s(K + \jj)^r\\
&&   \bip \bip_{i=0}^{\infty} \ii^s(Jt)^{2r+i} 
     \bip \bip_{i=1}^{r-1} \ii^s (K + \jj)^i  (Jt)^{2r-2i} 
     \bip \bip_{i=0}^{r-1} \ii^s H (K + \jj)^i (Jt)^{2r-2i-1}
\eeqn
for all $r \geq 1$ and for all $s \geq 0$.
Therefore
\beqn
&&   \f{(M^2 + \jj T)^r (\ii T)^s}{(M^2 + \jj T)^{r+1} (\ii T)^s}\\
&=&  \left[ \f{\ii^s}{ (K + \jj) \ii^s} 
     \bip   \f{\ii^s}{ H \ii^s} \right]
     \bip   
     \bip_{i=1}^{r-1} 
            \f{ (K + \jj)^i \ii^s}{(K + \jj)^{i+1} \ii^s}
     \bip \bip_{i=0}^{r-1} 
            \f{H (K + \jj)^i \ii^s}{H (K + \jj)^{i+1} \ii^s}\\
&&   \bip   \f{(K + \jj)^r \ii^s}{(K + \jj)^{r+1} \ii^s}
     \bip 
     \left[ 
            \f{J^{2r} \ii^s}{J^{2r} (K + \jj) \ii^s} 
     \bip   \f{J^{2r+1} \ii^s}{ J^{2r+1} H  \ii^s} \right]
     \bip \bip_{i=1}^{r-1} 
            \f{ (K + \jj)^i  J^{2r-2i}\ii^s}
              { (K + \jj)^{i+1}  J^{2r-2i}\ii^s}\\
&&   \bip \bip_{i=0}^{r-1} 
             \f{ H (K + \jj)^i J^{2r-2i-1}\ii^s}
               { H (K + \jj)^{i+1} J^{2r-2i-1}\ii^s}.
\eeqn
For large exponents the length of the  modules appearing in the above 
sum are polynomials. We are interested only in those modules whose length
 will contribute to the terms of total degree d in $r$ and $s$. 
In what follows, we will denote by $\cdots$ a function in $r$ and $s$
of total degree less than $d$. 
For large exponents 
$\ell ({\ii^s}/{ H  \ii^s})$,  
$\ell ({ J^{2r+1} \ii^s}/{ J^{2r+1} H \ii^s})$
and
$ \ell((K + \jj)^r \ii^s / (K + \jj)^{r+1} \ii^s)$
 are polynomials of total degree atmost $d-1$  [Lemma~\ref{kmv}]. 
Consider,
\beqn
&& \ell \left( \f{H (K + \jj)^i     J^{2r-2i-1} \ii^s}
                 {H (K + \jj)^{i+1} J^{2r-2i-1} \ii^s} \right) 
=  \ell \left( \f{  (K + \jj)^i     J^{2r-2i-1} \ii^s}
                   {(K + \jj)^{i+1} J^{2r-2i-1} \ii^s} \right) \\
&& 
+  \left[
   \ell \left( \f{  (K + \jj)^{i+1} J^{2r-2i-1} \ii^s}
                 {H (K + \jj)^{i+1} J^{2r-2i-1} \ii^s} \right)
-  \ell \left( \f{  (K + \jj)^{i}   J^{2r-2i-1} \ii^s}
                 {H (K + \jj)^{i}   J^{2r-2i-1} \ii^s} \right) \right].
\eeqn
By Lemma~\ref{kmv}, for large exponents the length of both the
modules appearing in the  square bracket are polynomials of total
degree at most $d-1$. Moreover, the coefficients of all the monomials
of highest degree appearing 
in both the polynomials are the same.  Hence,   their
difference is a 
polynomial of total degree at most $d-2$.  Thus
\beqn
      \ell \left( \f {H (K+{\jj})^i     J^{2r-2i-1} {\ii}^s} 
                     {H (K+{\jj})^{i+1} J^{2r-2i-1} {\ii}^s} \right) 
& = & \ell \left( \f {(K+{\jj})^i       J^{2r-2i-1} {\ii}^s}
                     {(K+{\jj})^{i+1}   J^{2r-2i-1} {\ii}^s} \right) 
                     + \cdots\\ 
& = & \ell \left( \f {(K+{\jj})^i       J^{2r-2i} {\ii}^s}
                     {(K+{\jj})^{i+1}   J^{2r-2i} {\ii}^s} \right) 
                      + \cdots . 
\eeqn
for large $r$ and $s$.
Similarly one can show 
\beqn
    \ell \left( \f{H (K + \jj)^i \ii^s}
                  {H (K + \jj)^{i+1} \ii^s} \right)
=   \ell \left( \f{  (K + \jj)^i \ii^s}
                  {  (K + \jj)^{i+1} \ii^s} \right) 
             + \cdots .
\eeqn
Thus considering the relevant terms we get
\beq
\label{eq1} \nonumber
&&     \ell \left( 
       \f{(M^2 + \jj T)^r\ii^sT}
         {(M^2 + \jj T)^{r+1}\ii^sT} \right) \\ 
&=&  2 \left[ \sum_{i=0}^{r-1}  
       \ell \left(
       \f{(K + \jj)^i \ii^s}
         {(K + \jj)^{i+1} \ii^s}\right)
   +  \sum_{i=0}^{r-1}  \ell \left(
       \f{(K + \jj)^i     J^{2r-2i} \ii^s}
         {(K + \jj)^{i+1} J^{2r-2i} \ii^s} \right) \right] + \cdots. 
\eeq
For large exponents  the terms which appear in the above sums are
polynomials of degree $d-1$. Without loss of generality, we can assume 
that they are polynomials for all exponents since the multiplicity
       formula will not be altered.  Consider 
\beq
\label{ns1} \nno
&&  \sum_{i=0}^{r-1}  \ell \left(
    \f{(K + \jj)^i     J^{2r-2i} \ii^s}
      {(K + \jj)^{i+1} J^{2r-2i} \ii^s} \right)\\ \nno 
&=& \sum_{i=0}^{r-1}
    \sum_{q_0+q_{1}+q=d-1} 
     e((K + \jj)^{[q_0+1]}|J^{[q_{1}]}|\ii^{[q]})
     {i+q_0 \choose q_0} 
     {2r-2i +q_{1} \choose q_{1}} 
     {s+q \choose q} + \cdots \\ \nno
&=& 2^{q_{1}}  
    \sum_{q_0+q_{1}+q=d-1} 
     e((K + \jj)^{[q_0+1]}|J^{[q_{1}]}|\ii^{[q]})
     {r + d-q \choose d-q}
       {s+q \choose q} + \cdots  \\
&& \hspace{3in}
    \mbox{[by  Lemma~\ref{prod} and  Lemma~\ref{bin}]}.
\eeq
Similarly, one can show that
\beq
\label{ns2}
    \sum_{i=0}^{r}
    \ell \left( \f{(K + \jj)^i \ii^s}{(K + \jj)^{i+1} \ii^s} \right)
&=& \sum_{q=0}^{d-1} e_{q}(K + \jj| \ii) 
           {r+d-q \choose d-q}  {s+q \choose q} + \cdots.
\eeq
Substitute (\ref{ns1}) and (\ref{ns2}) in   (\ref{eq1}).
By comparing the coefficient of  $r^{d-q} s^q$ ($0
\leq q \leq d-1$) in (\ref{eq1}) and (\ref{ns5}) we get the desired result.    
\qed

\vspace{.2in}

\bco
\label{extra2}
{\em \cite[Theorem~3.4]{katz-ver}}
Let $I$ be an ideal of positive height and let $J$ be an $\m$-primary
in $R$.  Let  $T = \b(I)$ and let  $M = (t^{-1}, J, It)$. Then 
\beqn
   e(M)
&=& \f{1}{2^{d}}
  \left[e (J^2 + I)
  + \sum_{q =0}^{d-1} 2^{q}
     e_{q}(J^2 + I | I) \right]
\eeqn
\eco
\pf Put $s=0$ in the above proof. 
\qed

\vspace{.2in}

\bco 
\label{l5}
Let $I$ be an $\m$-primary ideal and  let $J, I_{1}, \ldots, I_n \seq
\m$ be ideals of positive height in $R$. Let $\jj \seq I+J$ be any ideal
of $R$.  Put $T = \b(J)$, $M = (t^{-1}, I, Jt)$ and  $K = I^2
+ J$.  Then for all non-negative integers  $q, q_{1}, \ldots, q_n$
satisfying $q + q_{1} + \cdots + q_n = d-1$,
\beqn
&& e((M^2 + \jj T)^{[q_0+2]} | (I_{1}T)^{[q_{1}]} | \cdots | 
                       (I_nT)^{[q_n]})\\
&=& 2 \left[ e((K + \jj)^{[q_0+1]} | I_{1}^{[q_{1}]} | \cdots |
                                    I_n^{[q_n]})
 + \sum_{ k + l = q_0} 2^{l}
      e( (K + \jj)^{[k+1]}| J^{[l]} | I_{1}^{[q_{1}]}|
                 \cdots | I_n^{[q_n]}) \right].
\eeqn
\eco
\pf The proof follows by replacing 
$(I T)^{s}$ by $(I_{1}T)^{s_{1}} \ldots (I_{n}T)^{s_n}$ 
in (\ref{eq1}) of  Proposition~\ref{l4}  and by using arguments
      similar to those in Proposition~\ref{l4}.
\qed

\vspace{.2in}

\bl
\label{l7}
Let $I_{1}, \ldots, I_g \seq \m$ be
ideals of positive height in $R$. Put $I_{g+1} = (0)$.  Let $1 \leq j \leq g$.
 Then for all 
for all $j=1, \ldots, g$ and  for all non-negative integers
$q_0, q_{j+1}, \ldots,  q_g$ satisfying
$q_0 + q_{j+1} + \cdots + q_g = d-1-j$, 
\beqn
&&    e( L_{j}^{[q_0+j+1]} | I_{j+1} \b_{j}^{[q_{j+1}]} | 
                    \cdots | I_g \b_{j}^{[q_g]})   \\
&=&   2^j \sum_{t=0}^j 
     \sum_{\sta{\scriptstyle q + q_{1} + \cdots + q_t = q_0}
     {1 \leq i_{1} < \cdots < i_t \leq j}}
  2^{q_0  - q}
  e( L^{[q + 1]} | I_{i_{1}}^{[q_{1}]} | \cdots 
                  | I_{i_t}^{[q_{t}]} 
    | I_{j+1}^{[q_{j+1}]} | \cdots | I_g^{[q_g]}).
\eeqn
\el
\pf Notice that 
$
L_j = \n_j^2 + (I_{j+1} + \cdots + I_g)\b_j
$,
$
1 \leq j \leq g
$. 

We induct on $j$. Let $j = 1$. In  Corollary \ref{l5}, put  $n = g$, 
$T = \b_1$, 
$M = \n_1$, 
$J = I_1$, 
$I=\m$,
$J_1 = I_2 + \cdots + I_g$, 
and replace the set of ideals $\{I_1, \ldots, I_n \}$ by the set of
ideals  $\{I_2, \ldots, I_g\}$.
Also put $l = q_1$ and $k=q$. 

Suppose 
 $j>1$.  In  Corollary \ref{l5} put 
$T = \b_j$,  
$M = \n_j$ and   
$R = \c_{j-1} := (\b_{j-1})_{N_j}$
$ J = I_j$,
$I = \n_{j-1}$,
$J_1 = I_{j+1} + \cdots + I_g$
and replace the set of ideals $\{I_1, \ldots, I_n \}$ 
by the set of ideals $\{I_{j+1}, \ldots, I_g\}$.
Also  put  $l = q_j$ and $ k=q$. Then 

\beq
\label{m3} \nno
&&  e(L_j^{[q_0+j+1]} | I_{j+1}  \b_{j}^{[q_{j+1}]}| \cdots
                      | I_{g} \b_{j}^{[q_g]})\\ \nno
&=&  2 \left[ 
    e(L_{j-1}^{[q_0+j]}| I_{j+1} \c_{j-1}^{[q_{j+1}]} | \cdots | 
                               I_g \c_{j-1}^{[q_g]})  \right.\\ \nno
&&  \left. 
  + \sum_{q + q_{j} = q_0 + j-1} 2^{q_{j}}
    e( L_{j-1}^{[q+1]} |  I_{j} \c_{j-1}^{[q_{j}]} | 
                          I_{j+1} \c_{j-1}^{[q_{j+1}]}
              | \cdots | I_n \c_{j-1}^{[q_n]}) \right] 
%\\ \nno
 \hspace{.5in}  
\mbox{[by Corollary~\ref{l5}]} \\ \nno
&=&  2 \left[ 
    e(L_{j-1}^{[q_0 + j]} | I_{j+1} \b_{j-1}^{[q_{j+1}]} 
                 | \cdots | I_g \b_{j-1}^{[q_g]})  \right.\\  \nno
&&  \left. 
  + \sum_{q + q_{j} = q_0 +j-1} 2^{q_{j}}
    e( L_{j-1}^{[q + 1]} | I_{j} \b_{j-1}^{[q_{j}]} | I_{j+1} \b_{j-1}^{[q_{j+1}]}
            | \cdots | I_g \b_{j-1}^{[q_g]}) \right] 
%\\
 \hspace{.5in}  
\mbox{[by Lemma~\ref{dim2}]}.
\eeq
By induction hypothesis, each term in the above bracket can be
expressed as a sum of mixed multiplicities of ideals in the ring
$R$. Combining these terms in a nice way we get the desired result.
\qed

\vspace{.2in}

\noin
{\bf Proof of Theorem~\ref{mult}}
Since  $\n({\i})$ is a maximal ideal in $\b(\i)$, 
\beqn
  \ell \left( \f{\n({\i})^r \b(\i)_{\n({\i})}}{\n({\i})^{r+1} 
                            \b(\i)_{\n({\i})}} \right)
= \ell \left( \f{\n({\i})^r}{\n({\i})^{r+1}} \right).
\eeqn
By Lemma~\ref{dim1},  $\dim\b_g = d+g$. Hence
$ e(\n({\i})^2)
= 2^{d+g}e(\n({\i}))$.
Also $e(L_g^{[d+g]}) = e(L_g) = e(\n({\i})^2)$. 
Put $j = g$ in Lemma~\ref{l7}. This completes the proof of  the
theorem. 
\qed
\vspace{.2in}

\bco
\label{multone}
Let $(R,\m)$ be a local ring of dimension one. 
Let 
$I_{1}, \ldots, I_{g}$ be ideals of positive height.   Put $L = \m^2 + I_{1}
+ \cdots + I_{g}$. Then 
\beqn
   e(\b(\i)_{\n(\i)}) = 2^{g-1} e(L).
\eeqn
\eco
\pf Put $d=1$ in Theorem~\ref{mult}. Then
\beqn
 e(\b(\i)_{\n(\i)})
= \f{1}{2} \left[ 
  \sum_{t=0}^{g} \sum_{1 \leq i_{1} < \cdots < i_t \leq g} e(L) \right]
= \f{e(L)}{2} \left[ \sum_{t=0}^{g}  {g \choose t}      \right]
= 2^{g-1} e(L).
\eeqn
\qed

\vspace{.2in}

\bco
\label{multtwo}
Let $(R,\m)$ be a local ring of dimension two. 
Let 
$I_{1}, \ldots, I_{g}$ be ideals of positive height. Put 
$L = \m^2 + I_{1} + \ldots + I_{g}$.  Then 
\beqn
   e(\b(\i)_{\n(\i)}) 
= 2^{g-2} \left[ e(L) 
+ \sum_{j=1}^{g}e_{1}(L|I_{j}) \right].
\eeqn
\eco
\pf Put $d=2$ in  Theorem~\ref{mult}. Then
\beqn
      e(\b(\i)_{\n(\i)}) 
&=&  \f{1}{4} \left[ 
     \sum_{t = 0}^{g} 
     \sum_{\sta{\scriptstyle q + q_{1} + \cdots + q_t = 1}
          {1 \leq i_{1} <  \cdots < i_t \leq g}}
    2^{q_{_{1}} +  \cdots +  q_{t}}
    e({L^{[q+1]} | I_{i_{1}}^{[q_{1}]} | 
          \cdots | I_{i_t}^{[q_{t}]}}) \right]  \\
&=& \f{1}{4} \left[ 
    \sum_{t=0}^{g} {g \choose t} e(L) 
 +   2 \sum_{t=1}^{g-1} {g-1 \choose t-1} 
     \left[
     \sum_{j=1}^{g} e_{1}(L|I_{j}) \right] \right] \\
&=&  2^{g-2} \left[ e(L) + \sum_{j=1}^{g} e_{1}(L|I_{j}) \right].
\eeqn
\qed
\vspace{.2in}

We exhibit an  interesting   relationship
between the multiplicity formula of  the Rees  algebra 
$\r({\i})_{\mm(\i)}$ and that  of  the extended Rees algebra
$\b({\i})_{{\n}({\i})}$.

\vspace{.2in}

\brm
{\em 
Let  $I_{1}, \ldots, I_{g} \seq \m^2$ be  ideals of positive height in
$(R,\m)$.  Then 
\beqn
 e(\b({\i})_{{\n}({\i})})
= e(R) + 
   \sum_{t=1}^{g} 
  \sum_{1 \leq i_{1} <  \cdots < i_t \leq g}
   e(\r( I_{i_{1}},  \ldots, I_{i_t})_{\mm(I_{i_{1}},  \ldots, I_{i_t})}). 
\eeqn}
\erm

\vspace{.2in}

\brm
{\em  The multiplicity formula characterizes the
minimal multiplicity of $\b({\i})_{{\n}({\i})}$ (see \cite{clare}). }
\erm

\vspace{.2in}

\single

\noin
{\bf Acknowledgments:} {\small   The author is very grateful to 
J.~K.~Verma for his valuable suggestions and useful
conversations. The author also wishes to thank the National Board for
Higher Mathematics, D.~A.~E., India,  for financial support and the
Indian Institute of Technology, Bombay,  where the main work of 
this paper was carried out . The author thanks the referee for the
suggestions.}

\single

\end{document}